\documentclass[12pt, oneside]{article}   	
\usepackage{geometry}  
\usepackage{graphicx}				
\usepackage{amssymb}
\newcommand{\be}{\begin{equation}}
\newcommand{\ee}{\end{equation}}

\newcommand{\bftheta}{\mbox{\boldmath $\theta$}}
\newcommand{\Var}{\mbox{Var}}
\newcommand{\VarCov}{\mbox{Var-Cov}}
\title{Statistical Meaning of Mean Functions}
\author{Abram M. Kagan\\
Department  of Mathematics, University of Maryland\\
 College Park, MD 20742, USA\\
 Paul J. Smith\\
 Department  of Mathematics, University of Maryland\\
 College Park, MD 20742, USA}
\begin{document}
\maketitle
\begin{abstract}
\noindent
The basic properties of the Fisher information allow to reveal the statistical meaning of classical inequalities between mean functions.
The properties applied to scale mixtures of Gaussian distributions lead to a new mean function of purely statistical origin, unrelated to the classical arithmetic, geometric, and harmonic means. We call it the informational mean and show that when the arguments of the mean functions are Hermitian positive definite matrices, not necessarily commuting, the informational mean lies between the arithmetic and harmonic means, playing, in a sense, the role of the geometric mean that cannot be correctly defined in case of non-commuting matrices.\\
Surprisingly the monotonicity and additivity properties of the Fisher information lead to a new generalization of the classical inequality between the arithmetic and harmonic means.
\end{abstract}

\section{Introduction.}

Fisher information is a fundamental concept in statistics because it quantifies the efficiency of point estimators in finite samples and the asymptotic behavior of maximum estimators.  The importance of Fisher information is derived from two properties:
\begin{itemize}
\item {\it Monotonicity}: The Fisher information in a statistic (a reduction of a set of data) is never greater than the information in the complete data set.

\item {\it Additivity}: The total Fisher information in a set of independent observations is the sum of the Fisher informations of each of its components.
\end{itemize}

In this article we apply Fisher information to develop analytic inequalities involving both scalars and matrices. The monotonicity and additivity of Fisher information are key tools in deriving or reproving analytical inequalities, as shown below. Our general approach is to formulate a probability model, specialize it to Gaussian distributions, and use information-theoretic properties of the model to derive inequalities based on statistical principles.
 
In Kagan and Smith (2001) we used Fisher information to create statistical proofs of the monotonicity and convexity of the matrix function ${\bf A}^{-1}$ for Hermitian matrices.  That is, 
$$  {\bf A} \ge {\bf B} \Rightarrow {\bf B}^{-1} \ge {\bf A}^{-1}  $$ 
and, given weights $w_1,\ldots,w_n$ such that $w_j \ge 0$ and $\sum w_j =1$,
$$
  (w_1 {\bf A}_1 + \cdots + w_n {\bf A}_n)^{-1} \le 
               w_1 {\bf A}_1^{-1} + \cdots + w_n {\bf A}_n^{-1}.
$$
Here and throughout the paper, for any pair of Hermitian matrices, ${\bf A} \ge {\bf B}$ means ${\bf A} - {\bf B}$ is nonnegative definite.  Similarly the matrix function ${\bf A}^2$ is shown to be convex using statistical methods.

The convexity result above was extended to a notion of matrix-weighted averages in Kagan and Smith (1999).  The scalar weights in $w_1 {\bf A}_1 + \cdots + w_n {\bf A}_n$ are replaced by matrix weights as follows: 
$$
  {\bf B}_1^T {\bf A}_1{\bf B}_1 + \cdots + {\bf B}_n^T{\bf A}_n {\bf B}_n
$$
where ${\bf B}_1^T {\bf B}_1 + \cdots + {\bf B}_n^T {\bf B}_n = {\bf I}$.  It was shown that ${\bf A}^2$ and ${\bf A}^{-1}$ are {\it hyperconvex} functions, meaning that 
$$
  ({\bf B}_1^T {\bf A}_1{\bf B}_1 + \cdots + {\bf B}_n^T{\bf A}_n {\bf B}_n)^2
    \le 
  {\bf B}_1^T {\bf A}_1^2{\bf B}_1 + \cdots + {\bf B}_n^T{\bf A}_n^2 {\bf B}_n 
$$
   and
$$
  ({\bf B}_1^T {\bf A}_1{\bf B}_1 + \cdots + {\bf B}_n^T{\bf A}_n {\bf B}_n)^{-1}
    \le
  {\bf B}_1^T {\bf A}_1^{-1}{\bf B}_1 + \cdots 
                  + {\bf B}_n^T{\bf A}_n^{-1} {\bf B}_n. 
$$
As before, these results were derived by making use of the properties of Fisher information.
  
Our work is similar to the use of properties of entropy and related informational quantities to derive and extend classical inequalities.  See Dembo, Cover and Thomas (1991) for an exposition of that work.

\section{Properties of Fisher Information.}

Basic results concerning Fisher information are given in standard textbooks on mathematical statistics, for example Rao (1971) or Bickel and Doksum (2015).
Let $\bf X$ be a random vector with density $p({\bf x};\bftheta)$ depending on a parameter $\bftheta$. We assume the score function
$$
      J({\bf x};\bftheta) = (\partial/\partial\bftheta)\log p({\bf x};\bftheta)
$$ 
is well defined.  Then ${\bf I}_{{\bf X}}(\bftheta)$, the Fisher information on $\bftheta$ contained in $\bf X$, is defined as 
$$
   {\bf I}_{X}(\bftheta) = \VarCov[[J(X;\bftheta)] = E_{\bftheta}[J(X;\bftheta)J(X;\bftheta)^T].
$$ 
Under further regularity conditions, 
$$
   {\bf I}_{\bf X}(\bftheta) = E_{\bftheta}\left[-\frac{\partial}{\partial\bftheta} 
                                            \frac{\partial}{\partial\bftheta^T}\log p(x;\theta)\right].
$$  
The fundamental information inequality (or Cram\'{e}r-Rao inequality) states that if ${\bf T}$ is an unbiased estimator of $\bftheta$, then
$$
    \Var_{\bftheta}[{\bf T}] \ge {\bf I}_X(\theta)^{-1}.
$$
(If {\bf A} and {\bf B} are Hermitian matrices, the notation ${\bf A} \ge {\bf B}$ means that ${\bf A}- {\bf B}$ is nonnegative definite.)

When $\bftheta$ is a location parameter, $\bf X$ has density $p({\bf x}-\bftheta)$.  The Fisher information on a location parameter becomes 
$$ 
   {\bf I}_{\bf X} =\int(\partial\log p({\bf x})/\partial {\bf x})(\partial\log p({\bf x})/\partial {\bf x}^T) p(x)dx.
$$ 
Plainly, ${\bf I}_{\bf X}(\bftheta)={\bf I}_{\bf X}$ is constant in $\bftheta$. (The notation ${\bf I}_{\bf X}$  by default denotes the information on a location parameter $\bftheta$ throughout this paper.)

If ${\bf X}_\sigma$ is distributed as $\sigma {\bf X}$, the density of ${\bf X}_\sigma$ is $(1/\sigma)p(({\bf x}-\bftheta)/\sigma)$ and plainly ${\bf I}_{{\bf X}_\sigma} = {\bf I}_{\bf X} /\sigma^2$. 

For a scalar Gaussian random variable $X \sim N(\theta,\sigma^2)$ one has $I_X =1/\sigma^2$, and for any $X$ with
$E[X]=\theta$ and $\Var(X)=\sigma^2$, $I_X\geq 1/\sigma^2$.  This is a consequence of the Cram\'{e}r-Rao inequality.

\section{Mixtures, Mean Functions and Inequalities.}

Consider an experiment consisting of observing a pair $(\Delta,X)$, where $\Delta$ is a discrete random variable with $P(\Delta =i) = w_i$ and the conditional distribution of $X$ given $\Delta =i$ is $N(\theta,\sigma^{2}_i),\:i=1,\ldots,n$.

The marginal distribution of $X$ is a scale mixture of Gaussian distributions $N(\theta, 
\sigma^{2}_1),\ldots, N(\theta, \sigma^{2}_n)$ with mixture parameter ${\bf w}=(w_1,\ldots, w_n)$. Its density is
\be
p(x-\theta)=w_1\varphi_{\sigma_1}(x-\theta)+\cdots + w_n\varphi_{\sigma_n}(x-\theta).
\ee
Here $\varphi(x)$ is the density of the standard normal $Z \sim N(0,1)$.  The variance $\sigma^2$ of $X$ with density (1) is
\be
\sigma^2 = w_1 \sigma_{1}^2 +\cdots + w_n \sigma_{n}^2.
\ee 
The Fisher information on $\theta$ contained in the pair $(\Delta, X)$ is
\be
I_{(\Delta, X)}= w_1/\sigma^{2}_1 +\cdots+ w_n/\sigma^{2}_n.
\ee
Monotonicity of the Fisher information (the information in whole data set is never less than in any part of it; in our case $X$ is a part of $(\Delta, X)$) implies 
\[
                    I_X \leq I_{\Delta, X}.
\]
For any $Y$ with $E[Y]=\theta$, $I_Y \geq 1/\Var(Y)$. Hence one gets a two-sided inequality for $I_X$ with density $p(x-\theta)$:
\be 
 \left[\sum_{1}^n w_{i}\sigma^{2}_i\right]^{-1} \leq I_X \leq \sum_{1}^n w_i/ \sigma^{2}_i.
\ee
Since $p(x-\theta)$ in (1) is completely determined by the weights $w_1,\ldots,w_n$ and variances
$\sigma^{2}_1,\ldots,\sigma^{2}_n$, so is $I_X$. On setting $a_1=1/\sigma^{2}_1,\ldots,a_n=1/\sigma^{2}_n$,
the inequality (4) takes the form
\be
 \left[\sum_{1}^n w_{i}/a_i\right] \leq I_X(a_1,\ldots,a_n;w_1,\ldots,w_n) \leq \sum_{1}^n w_ia_i.
\ee

Recall that a function $M(a_1,\ldots, a_n)$  is called a {\it mean function} if for all\\ $a_1\geq 0,\ldots, a_n \geq 0$:
\begin{enumerate}
\item[(i)] $\min(a_1,\ldots,a_n)\leq M(a_1,\ldots,a_n)\leq \max (a_1,\ldots, a_n)$,
\item[(ii)]for any $\lambda >0$, $M(\lambda a_1,\ldots, \lambda a_n)= \lambda M(a_1,\ldots,a_n)$.
\end{enumerate}
Classical examples of mean functions are the arithmetic, geometric and harmonic means.

From (5), $I_X (a_1,\ldots,a_n;w_1,\ldots,w_n)$  satisfies (i).  Furthermore, for any $\lambda >0$, $I_X (\lambda a_1,\ldots,\lambda a_n;w_1,\ldots,w_n)$ is the Fisher information in $X_\lambda$ with density
\[p_{\lambda}(x-\theta)=w_1\varphi_{\sigma_1 /\lambda} +\ldots + w_n\varphi_{\sigma_n /\lambda}=\sqrt{\lambda}p(\sqrt{\lambda}(x-\theta))\]
and due to the well known property of the Fisher information mentioned above,
\[I_X(\lambda a_1,I\ldots,\lambda a_n;w_1,\ldots,w_n)= \lambda I_X(a_1,\ldots,a_n;w_1,\ldots,w_n)\]
so that $I_X(a_1,\ldots,a_n;w_1,\ldots,w_n)$ satisfies (ii).  Thus, $I_X(a_1,\ldots,a_n;w_1,\ldots,w_n)$ is a mean function.  We suggest calling it the {\it informational mean}.  

Inequalities (4) and (5) have a statistical interpretation. Their right hand sides are the Fisher information on $\theta$ in the pair $(\Delta, X)$ with
\be
P(\Delta=i)=w_i, ~~ X|\{\Delta=i\} \sim N(\theta,a_i=1/\sigma^{2}_i),~~ i=1,\ldots,n.
\ee
The left hand sides are the Fisher information on $\theta$ in a  Gaussian $X\sim N(\theta, \sigma^2)$ with
$\sigma^2$ given by (2).

Turn now to the case when $a_1,\ldots, a_n$ are replaced with Hermitian positive definite matrices ${\bf A}_1,\dots, {\bf A}_n$. As is well known, the inequality between the arithmetic and harmonic means still holds:
\be 
[w_1 {\bf A}^{-1}_1 +\cdots +w_n {\bf A}^{-1}_n]^{-1} \leq w_1 {\bf A}_1+\cdots+ w_n {\bf A}_n.
\ee
The matrices are not assumed to commute so that their geometric mean is not defined.

Suppose that $\bf X$ is a $d$-dimensional random vector  with distribution given by a density $p({\bf x}-\bftheta)$,
where ${\bftheta}=(\theta_1,\ldots,\theta_d)$ is a $d$-dimensional parameter, the vector score,
$$
    {\bf J}({\bf X}-\bftheta)
     = (\partial\log{p}/\partial\theta_1, \ldots, 
                        \partial\log{p}/\partial\theta_d)^T,
$$
is well defined and $E_{\bftheta}\|{\bf J}({\bf X}-\bftheta)\|^2 <\infty$.  Then the $d \times d$ matrix $E_{\bftheta}({\bf JJ}^T)=I_{\bf X}(\bftheta)$  is called the matrix of Fisher information on $\bftheta$ contained in $\bf X$.  (The superscript $T$ denotes transposition.)

For any Gaussian ${\bf Y} \sim N_d(\bftheta, {\bf V})$ with mean vector ${\bftheta}$ and non-degenerate covariance matrix ${\bf V}, I_{\bf Y}= {\bf V}^{-1}$.  For any ${\bf Y}\sim p({\bf y}-\bftheta))$ with covariance matrix
${\bf V}$, the information matrix is evidently constant in ${\bf \theta}$ and $I_{\bf X} \geq {\bf V}^{-1}$.  (Here and throughout this paper, ${\bf A} \geq {\bf B}$ means that the matrix ${\bf A-B}$ is nonnegative definite.)

Let $(\Delta, {\bf X})$ be a pair of random elements  whose distribution is given by
\be
P(\Delta=i)=w_i, ~~ {\bf X}|\{\Delta=i\} \sim  N(\bftheta,{\bf V}_i), ~~ i=1, \ldots, n.
\ee
The marginal density $p({\bf x}-\bftheta)$ of  $\bf X$ is the mixture of the densities of $N_{s}({\bftheta}, {\bf V}_1)$, $\ldots$, $N_{s}({\bftheta}, {\bf V}_n)$ with a mixture parameter $w_1,\ldots,w_n$.  Similarly to (2), the covariance matrix ${\bf V})$ of $\bf X$ is
\be
            {\bf V} = w_1 {\bf V}_1 +\cdots + w_n {\bf V}_n
\ee 
and the matrix of Fisher information on $\bf \theta$ in the pair $(\Delta,{\bf X})$ is
\be 
           I_{\Delta,{\bf X}} = w_1 {\bf V}^{-1}_1 +\ldots + w_{n} {\bf V}^{-1}_n,
\ee
which is constant in $\bftheta$.

As in the case of a scalar valued $\theta$, when $\bftheta$ is vector valued, the matrix of Fisher information is monotone. In our case, $I_{\bf X}\leq I_{\Delta,{\bf X}}$.

On setting ${\bf A}_1 ={\bf V}^{-1}_1,\ldots, {\bf A}_n ={\bf V}^{-1}_n$, $I_{\bf X}$ becomes a function of ${\bf A}_1,\ldots,{\bf A}_n$ and the mixing probabilities $w_1,\ldots,w_n$.  Comparing it with $I_{\Delta,{\bf X}}$ on one side and with the matrix of Fisher information in a Gaussian
$Z\sim N_{s}({\bf \theta}, V)$ on the other leads to
\be
               (w_1 {\bf A}^{-1}_1 +\ldots +w_n {\bf A}^{-1}_n)^{-1}
                       \leq I_{\bf X}({\bf A}_1,\ldots, {\bf A}_n;w_1,\ldots,w_n)
                       \leq w_1 {\bf A}_1+\ldots+w_n {\bf A}_n
\ee
 We want to emphasize that the matrices ${\bf A}_1,\ldots, {\bf A}_n$ are not assumed to commute.
 
 As a function of ${\bf A}_1,\ldots, {\bf A}_n$, $I_{\bf X}$ satisfies the above condition (ii) and the following version of (i): if a matrix $\bar {\bf A}$ and a positive matrix $\underline{{\bf A}}$ are such that $\underline{{\bf A}}\leq A_i \leq \bar {\bf A},\:i=1,\ldots,n$, then $\underline{A}\leq I_{\bf X} \leq \bar A.$  The statistical interpretation of (11) is the same as that of (4) and (5).
 
 \section{An inequality for Fisher information in sums of random variables.}
 
 In the previous section, we considered the Fisher information in a scale mixture of Gaussian densities to obtain analytic inequalities of mean functions.  In this section we follow a different approach by examining the Fisher information on  weighted location parameters in an independent sample of $n$ observations.  The model is as follows.

 For independent $X_1,\ldots, X_n$ with finite Fisher information and
 $w_1>0,\ldots, w_n >0,\:w_1+\ldots+w_n=1$, set
 \be
 U_i=X_i + w^{\alpha}_i\theta ,~~ i=1,\ldots,n.
 \ee
 The information in $U_i$ on $\theta$ equals $I_{U_i}=w_{i}^{2\alpha}I_{X_i}.$  Observe that for any constant $c>0$, the information in $U_i$ equals that in $cI_{U_i} $.
 
 Multiplying both sides of (12) by $w _{i}^{\beta}$ with $\beta=1-\alpha$ and taking the sum of the results gives
 \[U=\sum_{1}^{n}w_{i}^{\beta}U_i =\sum_{1}^{n}w_{i}^{\beta}X_i +\theta\]
 whence
 \be
 I_U=I_{\sum_{1}^{n}w_{i}^{\beta} X_i}.
 \ee 
 The information about $\theta$ in the vector $(w_{1}^{\beta}U_1,\ldots,w_{n}^{\beta}U_n)$ with independent components is the same as in the vector $(U_1,\ldots,U_n)$.
 Due to monotonicity and additivity of the Fisher information,
 \be 
 I_U=I_{\sum_{1}^{n}w_{i}^{\beta}U_i}\leq \sum_{1}^{n}I_{U_i}
 \ee
 whence
 \be
 I_{\sum_{1}^{n}w_{i}^{\beta}X_i}\leq \sum_{1}^{n}w_{i}^{2\alpha}I_{X_i}
 \ee  
 for $\alpha +\beta =1$. For $n=2, \alpha=\beta =1/2$ this inequality is known (e.g., see Dembo, Cover \& Thomas 1999, Theorem 13).
 
 When the $X_i$ are independent Gaussian variables with variances $\sigma_{i}^2 =1/a_i,$ 
 the sum $\sum w_{i}^{\beta}X_i$ has a Gaussian distribution with variance
 $\sum w_{i}^{2\beta}/a_i$ and (15) takes the form 
 %\end{document}
 \be
 \sum_{1}^{n} w_{i}^{2\alpha}a_i\geq \frac{1}{\sum w_{i}^{2\beta}/a_i}.
 \ee
 for $\alpha, \beta $ subject to $\alpha+\beta=1$.
 
 Replacing $2\alpha, 2\beta$ with $\alpha, \beta$ subject to $\alpha+\beta=2$
 gives a generalization, in a sense, of the classical inequality between the arithmetic and harmonic means:
 \be 
  \sum_{1}^{n} w_{i}^{\alpha}a_i\geq \frac{1}{\sum w_{i}^{\beta}/a_i}
 \ee
 for $\alpha+\beta=2$.
 
 \section{General comments}
           
 The paper reveals statistical meaning of classical mean functions (see in this connection Rao (2000), Kagan and Smith (2001), Kagan (2003), Kagan and Rao (2003)) and
 introduces a new one of purely statistical origin, called the informational mean. It leads to a new inequality similar to the classical inequality between the arithmetic, geometric and harmonic means and holds when the arguments of the mean functions are Hermitian positive definite matrices, not necessarily commuting in which case the geometric mean cannot be defined.\\
 The material of the paper can be used as a part of the chapter on the Fisher information in graduate courses in Statistics.   
 
\begin{center}
                       {\bf REFERENCES}
\end{center}

\begin{enumerate}
\item Bickel, P.J. and Doksum, K.A. (2015), {\it Mathematical Statistics} (Vol. 1, 2nd ed.), Boca Raton: CRC Press.

\item Dembo, A., Cover, T.M., and Thomas, J.A. (1991), ``Information Theoretic Inequalities,''  {\it IEEE Trans. Information Theory} 37, 1501-1518.

\item Kagan, A. and Smith, P.J. (1999), ``A Stronger Version of Matrix Convexity as Applied to Functions of Hermitian Matrices,''  {\it J. Inequal. \& Appl.}, 3, 143-152.

\item Kagan, A. and Smith, P.J. (2001), ``Multivariate Normal Distributions, Fisher Information and Matrix Inequalities,'' {\it Int. J. Math Educ. Sci. Technol.}, 32, 91-96.

\item Kagan, A. (2003), ``Statistical Approach to Some Mathematical Problems'', {\it Austrian J. Statist.}, 32(1-2), 71-83.

\item Kagan, A. and Rao, C. R. (2003), ``Some Properties and Applications of the Efficient Fisher Score'', {\it J. Statist. Plann. Inference}, 116, 343-352.

\item Rao, C. R. (2000), ``Statistical Proofs of Some Matrix Inequalities'', {\it Linear Algebra Appl.}, 321, 307-320.

\item Rao, C.R. (1973), {\it Linear Statistical Inference and Its Applications}, Hoboken, NJ: J. Wiley \& Sons.
\end{enumerate}
\end{document}